\documentclass[12pt]{article}
\usepackage{amsfonts}
\usepackage{amssymb}
\usepackage{amscd}
\usepackage{verbatim}
\usepackage{latexsym}

\begin{document}

\title{Finiteness of Selmer Groups and Deformation Rings }

\author{Chandrashekhar Khare, Ravi Ramakrishna}

\date{}

\maketitle

\newtheorem{theorem}{Theorem}
\newtheorem{lemma}[theorem]{Lemma}
\newtheorem{prop}[theorem]{Proposition}
\newtheorem{fact}[theorem]{Fact}
\newtheorem{cor}[theorem]{Corollary}
\newtheorem{example}{Example}
\newtheorem{conjecture}[theorem]{Conjecture}
\newtheorem{definition}[theorem]{Definition}
\newtheorem{quest}{Question}
\newtheorem{memo}{Memo}
\newcommand{\rhobar}{\overline{\rho}}
\newcommand{\Sha}{{\rm III}}
\newcommand{\SelX}{H^1_{ {\cal N}_v}(G_X,Ad^0 (\rhobar))}
\newcommand{\dualX}{H^1_{ {\cal N}^{\perp}_v}(G_X,Ad^0 (\rhobar)^*)}
\newcommand{\SelS}{H^1_{ {\cal N}_v}(G_S,Ad^0 \rhobar)}
\newcommand{\SelSp}{H^1_{ {\cal N}_v}(G_{S\cup S'},Ad^0 \rhobar)}
\newcommand{\SelSps}{H^1_{ {\cal N}_v}(G_{S\cup S'\backslash \{s'\}},Ad^0 \rhobar)}
\newcommand{\SelSQT}{H^1_{ {\cal N}_v}(G_{S\cup Q \cup T},Ad^0 (\rhobar))}
\newcommand{\dualS}{H^1_{ {\cal N}^{\perp}_v}(G_S,Ad^0 (\rhobar)^*)}
\newcommand{\dualSp}{H^1_{ {\cal N}^{\perp}_v}(G_{S \cup S'},Ad^0 (\rhobar)^*)}
\newcommand{\dualSps}{H^1_{ {\cal N}^{\perp}_v}(G_{S \cup S'\backslash \{s'\}},Ad^0 (\rhobar)^*)}
\newcommand{\dualSQT}{H^1_{ {\cal N}^{\perp}_v}(G_{S\cup Q \cup T},
Ad^0 (\rhobar)^*)}
\newcommand{\SelSQTqj}{H^1_{ {\cal N}_v}(G_{S\cup Q \backslash \{q_j \} \cup T},Ad^0 (\rhobar))}
\newcommand{\dualSQTqj}{H^1_{ {\cal N}^{\perp}_v}(G_{S\cup Q \backslash \{q_j \}\cup T},Ad^0 (\rhobar)^*)}

%\tableofcontents

\section{Introduction}

Consider a continuous, absolutely
irreducible mod $p$ representation
$$\rhobar:G_{\mathbb{ Q}} \rightarrow GL_2({\bf k}),$$ with $\bf k$
a finite field of characteristic $p$, and
$G_{\mathbb{Q}}$ the absolute Galois group of the field
of rational numbers $\mathbb{Q}$ and $p$ an odd
prime.

Suppose now that $\rhobar$ is {\em modular}, that is $\rhobar$ arises
from a Hecke eigenform $f$. Then we know there exists a
representation $\rho_f:G_{\mathbb{Q}} \rightarrow GL_2( {\bf{\cal O}}_f)$
where ${\bf {\cal O}}_f$ is the ring of integers of some
finite extension of ${\mathbb{Q}}_p$ whose reduction mod a
uniformiser is $\rhobar$.
One can associate to $\rho_f$ a Selmer group that is defined
via Galois cohomology. The calculation of the order of this 
Selmer group played a key role in the proof of the Shimura-Taniyama-Weil
conjecture.

In this paper we do {\em not} assume  $\rhobar$ is modular. For
instance, if
$\rhobar$ is even, i.e. $det\rhobar(c)=1$ where $c \in G_{\mathbb{Q}}$
is a complex conjugation, then $\rhobar$ is necessarily
not modular. The point of this paper is to find,
for $\rhobar$ odd or even,  a deformation
${\rho}:G_{\mathbb{Q}}  \rightarrow GL_2({\bf {\cal O}})$
of $\rhobar$ where ${\bf {\cal O}}$ is again the ring of integers
of some finite extension of ${\mathbb{Q}}_p$  such that the Selmer
group associated to ${\rho}$ is finite.

We remark that in many cases 
deformations ${\rho}$ of $\rhobar$ had been
constructed in [R2], but finiteness of the associated Selmer
did not follow from that work. 

In the last section we show how to remove the 
Generalized Riemann Hypothesis from Theorem $2$ of
[R1] and construct irreducible potentially semistable
$2$ dimensional $p$-adic Galois 
representations that are deformations of $\rhobar$, when $\rhobar$ is
odd,
and are ramified at infinitely many primes.
This method also applies to show how many even   
mod $p$ Galois representations can be deformed to characteristic
zero representations ramified at infinitely many primes.

We now consider the specifics of our situation.
Once and for all we fix  $\rhobar$  satisfying the hypotheses
below.
These hypotheses imply others that we need but use only implicitly.
\begin{itemize}
      \item $\rhobar$ and $Ad^0(\rhobar)$ are absolutely irreducible
            Galois representations.
      \item The (prime to $p$) Artin conductor $N(\rhobar)$ of $\rhobar$ is 
            minimal amongst its twists.
      \item We assume, as we may do, that the Serre weight
            $k(\rhobar)$  of $\rhobar$ is between $2$ and $p+1$.
      \item If $\rhobar$ is even then for the 
            decomposition group $G_p$ above $p$  we assume that $\rhobar 
            \mid_{G_p}$ is not twist equivalent to
            $\left(\begin{array}{cc} \chi &0 \\0 & 1\end{array}\right)$
            or twist equivalent to the indecomposable
            representation $\left(\begin{array}{cc} \chi^{p-2} & *\\
            0& 1\end{array}\right)$ where $\chi$ is the mod $p$
            cyclotomic character.
      \item If $\rhobar$ is odd we assume
             $\rhobar \mid_{G_p}$ is not
             twist equivalent to the trivial representation or the 
             indecomposable
             unramified representation given by
             $\left(\begin{array}{cc} 1 & * \\0 & 1\end{array}\right)$.
      \item $p\geq 7$ and the order of the projective image
            of $\rhobar$ is a multiple of $p$.
\end{itemize}
If one excludes the second part of the
last condition, then 
$\rhobar$ is allowed to have image of order prime to $p$. In this 
case it is easy to see that there is a 
unique `Teichm\"{u}ller' deformation of $\rhobar$ to $W({\bf k})$
that has image isomorphic to that of $\rhobar$.

Let $S$ be the set of containing
$p$, $\infty$ and the primes at which
$\rhobar$ is ramified.
Note that 
the primes other than $p$ and $\infty$ at which
$Ad^0(\rhobar)$ and $\rhobar$ are ramified are the same
by our assumptions.
Consider the minimal deformation ring $R_S$ (also
called $R_{\phi}$ in many references)
of $\rhobar$ and the universal minimal 
deformation $\rho_S^{univ}:G_{\mathbb{Q}}
\rightarrow GL_2(R_S)$. The ring $R_S$, like all the rings that will
occur in the course of this paper, is a 
complete Noetherian local ring (CNL) with residue field ${\bf k}$
that is a $W({\bf k})$-algebra, and all maps between such rings
will be maps of $W({\bf k})$ algebras that induce the identity
on residue fields. 

The universal deformation has the
property that it is unramified at the primes outside $S$ and
satisfies the following conditions:
\begin{itemize}
     \item If $\ell \neq p$, and $\rhobar(I_{\ell})$ 
            has order $c_l$ with $p \not | c_l$ then $\rho^{univ}_S(I_l)$
            has the same order $c_{\ell}$.
     \item If $\ell \neq p$ 
           is such that $\rhobar(I_{\ell})$ has order
           $c_{\ell}=p$, then $\rho_S^{univ}(I_{\ell})$ has an unramified
           quotient of rank one. (Note that if $p | c_l$ then $c_l=p$).
     \item There is also, for $\bar{\rho}$ odd,
           a condition at $p$ described in the
           `Local at $p$ considerations'
           section of [R2]. This is the condition of ordinariness
           or flatness etc. For $\rhobar$ even there is {\em no}
           local condition
           imposed at $p$.
     \item We fix  determinants of our deformations.
           There is an integer $k$ with $2 \leq k \leq p$
           and a character $\zeta$ of $G_p$ such that 
           $\bar{\rho} |_{I_p}$ is 
           isomorphic to either $\psi^{k-1} \oplus \psi^{p(k-1)}$
           or
           $\left(\begin{array}{cc}
           \chi^{k-1} & * \\0 & 1\end{array}\right)$ where
           $\psi$ is a fundamental character of level $2$.
           We write the global character
           $det(\rhobar)= \alpha \chi^{k-1}$ and
           fix the determinants of deformations we consider to be 
           ${\tilde{\alpha}}{\varepsilon}^{k-1}$.

\end{itemize}

Then it is expected, at least for $\rhobar$ odd,
that $R_S$ is a finite, flat $W({\bf k})$-module:
the finiteness of $R_S[{1 \over p}]$ over $Frac(W({\bf k}))$
is predicted by the Fontaine-Mazur conjectures (cf. [FM], but also see
Conjecture \ref{flatness} below).

We need the following definition/lemma to state 
our main theorem.

\begin{definition}\label{mega}
  Let $S'$ be a finite set of primes disjoint from $S$, such that for 
  $s' \in S'$, $s'$ is not $\pm 1$ mod $p$ and $\rhobar({\rm
  Frob}_{s'})$ has eigenvalues with ratio $s'$. There is a CNL
  $W({\bf k})$-algebra $R_{S \cup S'}$
  that is the universal ring for the (representable) 
  deformation problem which to any
  CNL $W({\bf k})$-algebra $A$, with maximal ideal ${\sf m}_A$, 
  assigns the isomorphism classes 
  of continuous representations $G_{\mathbb{Q}} \rightarrow 
  GL_2(A)$ of fixed determinant as above 
  that are 
\begin{itemize} 
\item minimal at $S$ \item
   unramified outside $S \cup S'$, 
\item mod ${\sf m}_A$, the maximal ideal of $A$, {\bf equal}
  to $\rhobar$. \end{itemize} Denote by $\rho_{S \cup S'}^{univ}$ 
  the corresponding
  universal representation 
  $\rho_{S \cup S'}^{univ}:G_{\mathbb{Q}} \rightarrow
  GL_2(R_{S \cup S'})$. For any $\alpha \subset S'$ consider
  the deformation problem which to any
  CNL $W({\bf k})$-algebra $A$ assigns the isomorphism classes 
  of continuous representations $\rho:G_{\mathbb{Q}} \rightarrow 
  GL_2(A)$ of fixed determinant as above 
  that are \begin{itemize} \item minimal at $S$ \item
   unramified outside $S \cup S'$, 
\item mod ${\sf m}_A$, the maximal ideal of $A$, {\bf equal}
  to $\rhobar$ 
\item at primes $s' \in \alpha$, $\rho|_{G_{s'}}$ is up to
  twist of the form 
  $\left(\begin{array}{cc} \varepsilon &* \\0 & 1\end{array}\right)$.
\end{itemize} 
Then this deformation problem is representable (see
  lemma below) by
  a CNL $W({\bf k})$-algebra $R_{S \cup S'}^{\alpha-new}$, and we denote 
  by $\rho_{S \cup S'}^{\alpha-new}$ the corresponding
  universal representation
  $\rho_{S \cup S'}^{\alpha-new}:G_{\mathbb{Q}} \rightarrow
  GL_2(R_{S \cup S'}^{\alpha-new})$. There is a canonical 
  surjection $R_{S \cup S'}
  \rightarrow R_{S \cup S'}^{\alpha-new}$ for all $\alpha \subset S'$. 
  We denote the kernel of the map $R_{S \cup S'}
  \rightarrow R_{S \cup S'}^{S^{'}-new}$ by  
  $\wp_{S'}$.
\end{definition}

Below is the implicit representability claim
in Definition \ref{mega}. 

\begin{lemma}\label{newexists}
  The deformation problem where we require the properties
  of Definition \ref{mega} above is representable by
  a CNL $W({\bf k})$-algebra $R_{S \cup S'}^{\alpha-new}$, and we denote
  by $\rho_{S \cup S'}^{\alpha-new}$ the corresponding
  universal representation
  $\rho_{S \cup S'}^{\alpha-new}:G_{\mathbb{Q}} \rightarrow
  GL_2(R_{S \cup S'}^{\alpha-new})$. There is a surjection $R_{S \cup S'}
  \rightarrow R_{S \cup S'}^{S'-new}$ whose kernel we denote by
  $\wp_{S'}$.
\end{lemma}

\noindent {\bf Proof.}  We simply sketch the proof as it is quite
standard.
Using Mazur's work ([Ma2]) 
we know $R_{S \cup S^{'}}$ exists. Mazur also showed
for the prime $p$ (see [Ma1] for instance) 
that one could impose the condition 
of ordinariness at $p$.
Then one had a closed subfunctor of the original representable deformation
functor. Here we are imposing for each $s' \in \alpha$ a condition
almost identical to ordinariness. That $s'$ is not $1$ mod $p$ and
$\rhobar |_{G_{s'}}$ is of the form above is sufficient to guarantee that
a corresponding closed subfunctor of the original deformation functor exists 
(as these are deformation conditions in the sense of Sections 23 to 25
of [Ma1]). Letting $s'$ run through $\alpha$ and taking the intersection of 
these closed subfunctors, we are done.
 
\vspace{3mm}

The following theorem is the main result of this
paper (see Theorem \ref{main1} below):

\begin{theorem}\label{main}
      There is a finite set of primes
      $S'=\{s'_1,\dots, s'_n\}$
      not congruent
      to $\pm 1 \ mod \ p$, at which $\rhobar$ has eigenvalues with ratio
      $s'_i$ such that
      $R_{S \cup S'}^{S'-new} \simeq W({\bf k})$.
      We can choose the set
      of primes $S'$ so that the universal deformation
      $\rho_{S \cup S'}^{S'-new}:
      G_{\mathbb{Q}} \rightarrow GL_2(W({\bf k}))$
      corresponding to  $R_{S \cup S'}^{S'-new} \simeq W({\bf k})$ 
      is ramified at all the primes of $S'$.
\end{theorem}

\noindent{\bf Remarks:} 
\newline\noindent
$1$) The existence of $S'$ so that 
$R_{S \cup S'}^{S'-new} \simeq W({\bf k})$ follows 
directly from the methods of [R2]. 
However in [R2] it was not known that the representation
$\rho^{S^{'}-new}_{S \cup S'}$ was ramified at the primes of $S^{'}$.
For $s' \in S^{'}$, the corresponding versal local representations are  
indeed ramified, but it is not automatic that
$\rho_{S \cup S'}^{S'-new}$ is ramified  at all primes of $S'$.
Indeed, for some prime
$s' \in S^{'}$ we might have had that 
$\rho^{S^{'}-new}_{S \cup S'}|_{G_{s'}}$
was equal, up to twist, to 
$\left(\begin{array}{cc} \varepsilon & 0\\0 & 1\end{array}\right)$.
In the present work we obtain ramification at these primes by choosing
$S^{'}$ to be (possibly) larger than the ``auxiliary set'' of
[R2]. The main dichotomy that this paper addresses is
that, on the one hand the deformation
problem which assigns to $S'$ as above lifts of $\rhobar$ that are
{\bf ramified} at all primes in $S'$ is not representable (being ramified
at a prime is not a condition in the sense of Section 23 of [Ma1]), 
while on the other hand
it is not obvious that the universal representation 
$\rho_{S \cup S'}^{S'-new}$ {\it is ramified} at all primes
in $S'$.
\newline\noindent
$2$) For (many) odd $\rhobar$, the recent
important work of Taylor, cf.~[T1],
implies that $\rho_{S \cup S'}^{S'-new}$ is ramified at all primes
in $S'$ for {\bf all} $S'$ such that $R_{S \cup S'}^{S'-new} 
\simeq W({\bf k})$. Taylor's methods are partly geometric and 
very different from the methods of this paper.
The methods of this paper are purely Galois cohomological and work for 
$\rhobar$ odd and even uniformly, while Taylor's work does not 
address the case of $\rhobar$ even. In [T2] the results of [R2] are
extended to the case of $p=3,5$: then our methods do prove an analog of Theorem
\ref{main} for $p=5$. 
\vskip1em

Using Theorem 1 we can deduce easily (see Theorem \ref{finiteness1} below):

\begin{theorem}\label{finiteness}
 The tangent space $\wp_{S'}/\wp_{S'}^2$ is finite.
\end{theorem}

\noindent{\bf Remark:}
\newline\noindent
The tangent space above in Theorem \ref{main}
is called a Selmer group because of the
cohomological interpretation of its dual as in Proposition $1.2$ of 
[W] (see
also proof of
Theorem \ref{finiteness1}), hence the title of the paper. When $\rhobar$ is even
this Selmer group is not the one considered by Bloch and Kato, but it is the one
that is more pertinent when studying deformation rings.
In the case of $\rhobar$ even, in some particular cases,
deformations of $\rhobar$ to characteristic 0 were produced, prior to [R2], 
in [R4]. Using particular features of these lifts, 
the finiteness of the associated Selmer groups was
proven in [R4] and [B].  

\vspace{3mm}

We recall that an irreducible component of a ring $R$ is a quotient
$R/P$ for $P$ a minimal prime ideal.

\begin{cor}
  The prime ideal $\wp_{S'}$ is minimal. Thus the deformation
  ring $R_{S \cup S'}$ contains an 
  irreducible component isomorphic to $W({\bf k})$.
\end{cor}
We will also give an application of Theorem \ref{finiteness} to a conjecture of
Fontaine-Mazur (see Proposition \ref{fm} below).

We thank the anonymous referee for careful reading of
the original manuscript and many helpful suggestions.
In particular, the referee suggested a major simplification
of the results of section $3$.

\section{Auxiliary sets}

\begin{definition}\label{auxiliary}
  A finite set of primes $S'$ disjoint from $S$
  is said to be auxiliary if 
  for all primes $s' \in S'$, $s'$ is not $\pm 1$ mod $p$, $\rhobar({\rm
  Frob}_{s'})$ has eigenvalues with ratio $s'$, and  $R_{S \cup
  S'}^{S'-new}$ (defined in Definition \ref{mega}) is $\simeq W({\bf k})$. 
\end{definition} 

\begin{definition}\label{desire}
      If $\rho:G_{\mathbb{Q}} \rightarrow 
      GL_2(W({\bf k})/p^n)$ is a continuous
      representation that is a lift of $\rhobar$, 
      with determinant fixed as in the 
      introduction, then $\rho$ is special at $s'$, if $s'$ is 
not $\pm 1$ mod $p$ and  if $\rho|_{G_{s'}}$ can
      be conjugated to (we will also say for simplicity, is of the form)  
       $$\left(\matrix{\varepsilon&*\cr
                  0&1\cr}\right)$$ up to twist. The $*$ may be trivial.
\end{definition}

The $*$ can be genuinely nontrivial 
if and only if $\rho$ is ramified.

Recall from [R2] and [T2] that when studying a global deformation
problem, we prescribe a local deformation condition at each
of the primes at which ramification may occur. In particular,
for a local problem at $v$ we choose a class of deformations
of $\rhobar \mid_{G_v}$ to $W({\bf k})$ which we call ${\cal C}_v$.
For $v\in S$, ${\cal C}_v$ consists of deformations
of $\rhobar \mid_{G_v}$ to $W({\bf k})$ that 
are minimal in the sense of the
introduction.
Another way to think about this is that
${\cal C}_v$ is the set of $W({\bf k})$-valued
points on a certain smooth quotient of the unrestricted local 
deformation ring. 
The surjective map onto this smooth quotient induces
a surjective map on the (mod $p$) tangent spaces to these rings.
Since the dual of the tangent space of the unrestricted
local deformation ring is
$H^1(G_v,Ad^0(\rhobar))$, the smooth quotient  
gives rise to a subspace of $H^1(G_v,Ad^0(\rhobar))$ which
we denote
${\cal N}_v$. 

For each possible $\rhobar |_{G_v}$
one needs to compute  ${\cal C}_v$ and ${\cal N}_v$. In particular,
to find a characteristic zero deformation as in [R2] and [T2], 
one needs that 
$$dim_{\bf k} H^1(G_v,Ad^0(\rhobar))=dim_{\bf k}{\cal N}_v
+ dim_{\bf k} H^2(G_v,Ad^0(\rhobar)) +\delta$$
where $\delta =0$ except when $\rhobar$ is odd and $v=p$,
in which case $\delta=2$.
The calculations  of ${\cal N}_v$ and ${\cal C}_v$ are
laborious and not included here. See [R2].
Note, however, that the excluded situations of the introduction
are those for which the above equation does not hold.
(More precisely, at present we do not see how to choose 
${\cal C}_v$ and ${\cal N}_v$ appropriately in these cases).
For $\rhobar$ odd and not in an excluded case of the introduction,
${\cal C}_p$ consists of certain potentially semistable
representations. Locally, and therefore globally, these
representations are semistable up to a twist by a character of finite
order. For $\rhobar$ even the local deformation ring at $G_p$ is smooth
and ${\cal C}_p$ is taken to be the ($W({\bf k})$-valued
points of) the entire deformation
ring. Throughout this paper, for $s'$ a prime in an auxiliary
set, we will insist $\rhobar$ is special at ${s'}$
and 
${\cal C}_{s'}$ will consist of the special $W({\bf k})$-valued
points of the local at $s'$ deformation ring.
In this case elements of $H^1(G_{s'},Ad^0(\rhobar))$
are trivial on wild inertia, so we may consider them as functions
on tame inertia, which is topologically generated
by $\sigma_{s'}$ corresponding to Frobenius and $\tau_{s'}$
corresponding to (tame) inertia. 
${\cal N}_{s'}$ is spanned by 
the $1$-cohomology class given by
$$g(\sigma_{s'}) =\left(\begin{array}{cc}0&0\\0&0\end{array}\right)
\,\,\, g(\tau_{s'})=\left(\begin{array}{cc}0&1\\0&0\end{array}\right).$$
Also see the discussion after Lemma $2$ of [R3].

Recall from section $6$ of [R3]
that there is a subfield ${\bf \tilde{k}}$ of ${\bf k}$
which is the minimal field of definition
of the representations of $G_{\mathbb{Q}}$ on
$Ad^0(\rhobar)$ and $Ad^0(\rhobar)^*$. 
Denote the descents to ${\bf \tilde{k}}$ of $Ad^0(\rhobar)$ 
and $Ad^0(\rhobar)^*$
by $\widetilde{Ad}^0(\rhobar)$ and $\widetilde{Ad}^0(\rhobar)^*$ respectively.
 
\begin{lemma}\label{disjointness}
Let $\rhobar:G_{\mathbb{Q}} \rightarrow GL_2({\bf k})$ be given
with $S$ defined as usual. Assume $\rhobar$ satisfies the hypotheses
of the introduction and $X\supset S$ is a finite
set such that $X \backslash S$ consists of special primes.
For $n \geq 2$, let  $\rho_n : G_{X} \rightarrow
GL_2(W({\bf k})/p^{n})$ be a deformation of $\rhobar$
unramified outside $X$.
Let $\psi_i$ ($i=1,\cdots,m$) be finitely many 
${\bf \tilde{k}}$-linearly independent elements in 
$H^1(G_X,\widetilde{Ad}^0(\rhobar))$ 
and $\phi_j$($j=1,\cdots,r$) be finitely 
many ${\bf \tilde{k}}$-linearly independent elements
in 
$H^1(G_X,\widetilde{Ad}^0(\rhobar)^*)$.
Let $\mathbb{Q}(Ad^0(\rhobar))$ be the field
fixed by the kernel of the action of $G_{\mathbb{Q}}$
on $Ad^0(\rhobar)$. Let ${\bf K}=\mathbb{Q}(Ad^0(\rhobar),\mu_p)$
be the field obtained by adjoining the $p$th roots of unity to
$\mathbb{Q}(Ad^0(\rhobar))$.
We denote by ${\bf K}_{\psi_i}$ and ${\bf K}_{\phi_j}$
the fixed fields of the kernels of the restrictions
of $\psi_i, \phi_j$ to $G_{{\bf K}}$, 
the absolute Galois group of ${\bf K}$.
Let ${\bf P}_n$
be the fixed field of the kernel of the restriction of 
the projectivisation of $\rho_n$ to $G_{{\bf K}}$.
Then each of the fields ${\bf K}_{\psi_i}$, ${\bf K}_{\phi_j}$
${\bf P}_n$ and ${\bf K}(\mu_{p^n})$ is
linearly disjoint over ${\bf K}$ with the compositum of the others.
Let $I$ be a subset of $\{1,\cdots,m\}$ and $J$ a subset of $\{1,\cdots,r\}$.
Then there exists a prime $w \not \in X$ such that
\newline
\noindent
1) $\rho_{n-1}$, the mod $p^{n-1}$ reduction
    of $\rho_n$,  is special at $w$ but $\rho_n$ is {\em not}
   special at $w$. 
\newline
\noindent
2) $\psi_i |_{G_w} \neq 0$  for $i \in I$ and 
 $\psi_i |_{G_w} = 0$ for $i \in \{1,\cdots,m\} \backslash I$.
\newline
\noindent
3) $\phi_j|_{G_w} \neq 0$ for $j \in J$ and 
$\phi_j|_{G_w} =0$ for $j \in \{1,\cdots,r\} \backslash J$.
\end{lemma}
{\bf Proof.} 
(We actually require in $3$) above that $\phi_j
\not \in {\cal N}^{\perp}_w$, but since $\phi_j$ is unramified
at $w$, it is an exercise to see
this is equivalent to $\phi_j |_{G_w} \neq 0$).

The proof is quite standard 
and is essentially contained in proof of Lemma 1.2
of [T2].
The main differences are 
\begin{itemize}
   \item In [T2] $w$ is chosen only so that
         $\psi_i |_{G_w} \neq 0$ and 
         $\phi_j|_{G_w} \neq 0$.
   \item $\rho_n$ and condition $1$) do not figure there.
\end{itemize}
We take care of these 
in the proof we sketch below (as in Lemma 2 of [R1]).
Observe the field fixed by  kernel of the action of $G_{\mathbb{Q}}$
on  $Ad^0(\rhobar)^*$ is contained in ${\bf K}$. 
The following facts are in section $7$ of [R3].
See also Section 1 of [T2].

\begin{itemize}
   \item ${\bf K}_{\phi_j}$ and ${\bf K}_{\psi_i}$ are Galois over
          $\mathbb{Q}$ and linearly disjoint over ${\bf K}$
   \item $Gal({\bf K}_{\psi_i}/\mathbb{Q}) $ injects into 
         $Ad^0(\rhobar)$ and the map is $Gal({\bf K}/{\mathbb{Q}})$
         equivariant. Also the short exact sequence
         $$ 1\rightarrow  Gal({\bf K}_{\psi_i}/{\bf K}) \rightarrow
         Gal({\bf K}_{\psi_i}/\mathbb{Q}) \rightarrow 
          Gal({\bf K}/{\mathbb{Q}}) \rightarrow 1$$ splits.
   \item $Gal({\bf K}_{\phi_j}/{\bf K}) $ injects into
         $Ad^0(\rhobar)^*$ and the map is $Gal({\bf K}/{\mathbb{Q}})$
         equivariant. Also the short exact sequence
         $$ 1\rightarrow  Gal({\bf K}_{\phi_j}/{\bf K}) \rightarrow
         Gal({\bf K}_{\phi_j}/{\bf K}) \rightarrow
          Gal({\bf K}/{\mathbb{Q}}) \rightarrow 1$$ splits.
\end{itemize}

The $Gal({\bf K}/\mathbb{Q})$
equivariant isomorphisms $Gal({\bf K}_{\phi}/{\bf K}) \rightarrow
Ad^0(\rhobar)^*$ and $Gal({\bf K}_{\psi}/{\bf K}) \rightarrow
Ad^0(\rhobar)$ give ${\bf \tilde{k}}$ structures to
$Gal({\bf K}_{\psi_i}/{\bf K})$ and  
$Gal({\bf K}_{\phi_j}/{\bf K})$ 
making these maps isomorphisms of 
${\bf \tilde{k}}[Gal({\bf K}/\mathbb{Q})]$-modules.
We have $Gal({\bf K}_{\psi}/{\bf K}) \simeq
\widetilde{Ad}^0(\rhobar)$ and  $Gal({\bf K}_{\phi}/{\bf K}) \simeq
\widetilde{Ad}^0(\rhobar)^*$ as ${\tilde {\bf k}}[Gal({\bf K}/
\mathbb{Q}]$-modules.
Let $\widetilde{\rho_t}:G_{\mathbb{Q}}\rightarrow PGL_2(W({\bf k})/p^t)$
for $t=1,\cdots,n$ be the mod $p^t$ reduction
of the projectivisation $\widetilde{\rho_n}$
of $\rho_n$. Then we claim that the image of
$\widetilde{\rho_t}$ contains $PSL_2(W({\tilde {\bf k}}))/p^t$
This is true for $t=1$ by assumptions on $\rhobar$ and then
follows inductively for all $t$ by the proof of Lemma 3
on IV-23 of Serre's book [S] (here we use $p>3$).

Thus we have an exact sequence
$$1 \rightarrow \widetilde{Ad}^0(\rhobar) \rightarrow
Image(\widetilde{\rho_t})
\rightarrow Image(\widetilde{\rho_{t-1}}) \rightarrow 1.$$
It is easy to see, using the assumption that $p>3$, that 
this is a non-split exact sequence for all $t \geq 2$.

From this and the argument at the end of proof of Lemma 1.2 of [T2]
it follows easily that each of ${\bf K}_{\phi_j}$,
${\bf K}_{\psi_i}$, ${\bf P}_n$ and ${\bf K}(\mu_{p^n})$
is linearly disjoint over ${\bf K}$ with the compositum
of the others.
After this getting a $w$ that satisfies the conditions of the lemma
is an application of Chebotarev's theorem using
the (strong) linear disjointness of ${\bf K}_{\phi_j}$, ${\bf
K}_{\psi_i}$,
${\bf P}_n$ and ${\bf K}(\mu_{p^n})$
over ${\bf K}$.

\vspace{3mm}

\begin{definition} For a finite set $X$ containing $S$ we define
the {\em Selmer group}  
$\SelX$
to be the kernel of the map
$$H^1(G_{X },Ad^0(\rhobar)) \rightarrow \oplus_{v \in X}
H^1(G_v,Ad^0(\rhobar))/{\cal N}_v$$
and the {\em  dual Selmer group} $\dualX$ to be the kernel of the map
$$H^1(G_{X },Ad^0(\rhobar)^*) \rightarrow \oplus_{v \in X}
H^1(G_v,Ad^0(\rhobar)^*)/{\cal N}^{\perp}_v$$
where ${\cal N}^{\perp}_v \subset
H^1(G_v,Ad^0(\rhobar)^*)$ is the annihilator
of ${\cal N}_v \subset H^1(G_v,Ad^0(\rhobar))$ via local
duality.
\end{definition}

\begin{fact} For our given $\rhobar$ an auxiliary set exists.
\end{fact}
{\bf Proof.} We only give a brief sketch here. We refer the reader
to [R2] and [T2] for details.
There is a Galois cohomological
criterion (realisable!) for auxiliary sets
to exist, namely $\dualSp=0$.
Recall for any $X \supseteq S$ we get, from the Poitou-Tate
exact sequence,  the exact sequence
$$ H^1(G_X, Ad^0(\rhobar)) \rightarrow
\oplus_{v \in X} H^1(G_v,Ad^0(\rhobar))/{\cal N}_v
\rightarrow {\dualX}^{\vee}$$
$$\rightarrow H^2(G_X, Ad^0(\rhobar))\rightarrow
\oplus_{v \in X} H^2(G_v,Ad^0(\rhobar)).$$
If $\dualX$ is trivial then one can construct
a deformation $\rho$ of $\rhobar$ to $W({\bf k})$ such that
$\rho |_{G_v} \in {\cal C}_v$ for all $v \in X$
by deforming from mod $p^n$ to mod $p^{n+1}$.
One does this by observing the obstruction to deforming
lies in $H^2(G_X, Ad^0(\rhobar))$ and since
$\dualX=0$ this can be analyzed locally. The local
obstructions can be made to vanish by `moving'
from our given mod $p^n$ deformation to an unobstructed
mod $p^n$ deformation. That this is possible follows
from the surjectivity of
$$H^1(G_X, Ad^0(\rhobar)) \rightarrow
\oplus_{v \in X} H^1(G_v,Ad^0(\rhobar))/{\cal N}_v.$$
In fact the set of all such deformations ends up being
a power series ring over $W({\bf k})$ in
$dim\left(\SelX\right)$ variables.

\vspace{3mm}
\noindent
{\bf Remarks:} 
\newline\noindent
$1$) As all primes $s'$ in our 
auxiliary sets are special for $\rhobar$ we will
always have that $H^0(G_{s'}, Ad^0(\rhobar))$ and
$H^1(G_{s'}, Ad^0(\rhobar))/{\cal N}_{s'}$ are one dimensional.
See the Section $3$ of [R3].       
\newline\noindent
$2$) By Proposition $1.6$ of [W]
$\dualX$ and $\SelX$ have the same dimension.
Thus when we have annihilated the  Selmer group,
we have also annihilated the dual Selmer group
and the power series ring referred to in the proof above
is just $W({\bf k})$.
\newline\noindent
$3$)
It follows from [R2] and [T2]
  that any set $S'$, such that the primes
  in $S'$ are not $\pm 1$ mod $p$, and are special for
  $\rhobar$, and such that $R_{S \cup S'}^{S'-new}=W({\bf k})$,
  has cardinality at least $n_{\rhobar}=
  dim \,H^1_{{\cal N}^{\perp}_v}(G_S,Ad^0(\rhobar)^*)$.
We call these $S'$ {\bf minimal auxiliary sets}.

\section{Ramified auxiliary sets}

We begin with a definition.

\begin{definition}\label{ramified}
  An auxiliary set $S'$ such that the corresponding representation
  $\rho^{S^{'}-new}_{S \cup S'}$ 
  is ramified at all primes of $S'$ is said to be a
  ramified auxiliary set.
\end{definition}

The purpose of this section is to prove that {\bf ramified
auxiliary sets} do exist. 
We state the main theorem of this section which
is essentially a restatement of Theorem
\ref{main}.

 \begin{theorem}\label{main1}
   For our given $\rhobar$ there 
   are infinitely many ramified auxiliary sets. 
\end{theorem}

Consider the unique deformation 
$\rho^{{S'}-new}_{S \cup S'}:G_{\mathbb{Q}}
\rightarrow GL_2(W({\bf k}))$. If it is ramified
at all primes in $S'$ we are done.
If not, 
we write $S'= S'_{ram} \cup S'_{un}$ as the disjoint
union of primes that are ramified in
$\rho^{{S'}-new}_{S \cup S'}$ and primes that are
unramified in $\rho^{{S'}-new}_{S \cup S'}$.
We assume $S'_{un}$ is nonempty and if we can
discard any primes from $S'_{un} \subset S'$ and
the set remains auxiliary we do so.
Choose $s' \in S'_{un}$.
We will replace $s'$ by two primes $s_1$ and $s_2$
such that $S''=S' \cup \{s_1,s_2\}\backslash \{s'\}$
is an auxiliary set and
$\rho^{S''-new}_
{S \cup S''}$
is ramified at $\{s_1,s_2\} \cup S'_{ram}$.
Thus we replace $S'$ by an auxiliary
set that has one more element, but the set of
unramified primes for the new representation
is at least one smaller than for $S'$.
Repeating this process leads to a ramified auxiliary set.

Let $n$ be an integer such that
$\rho^{S'-new}_{S \cup S',n-1}$,
the mod $p^{n-1}$ reduction of
$\rho^{S'-new}_{S \cup S'}$,
is ramified at
all primes of $S'_{ram}$.
We will
consider the set $S' \backslash \{s'\} $.
The new deformation we construct will, mod $p^{n-1}$,
be equal to $\rho^{{S'}-new}_{S \cup S',n-1}$
and hence ramified at all primes of $S'_{ram}$.

\begin{prop} \label{h1surj}We assume $S \cup S' 
\backslash \{s'\}$ is not auxiliary. Let $\psi$
span the Selmer group for 
$S \cup S' \backslash \{s'\}$.
There exist primes  
$s_1$ and $s_2$ such that 
\newline\noindent
1) $\rho^{ {S'}-new}_{S \cup S',n-1}$ is special
at $s_1$
but $\rho^{ {S'}-new}_{S \cup S',n}$ is {\em not}
special
at $s_1$.
\newline\noindent
2) $\rho^{ {S'}-new}_{S \cup S',n-1}$ is special
at  $s_2$.
\newline\noindent
3)  The map 
$H^1(G_{S \cup S' \cup \{s_1\}  \backslash \{s'\} },Ad^0(\rhobar))
  \rightarrow \oplus_{v \in S \cup S' \backslash \{s'\}}
  H^1(G_v,Ad^0(\rhobar))/{\cal N}_v$ has one dimensional
kernel spanned by $\psi$
\newline\noindent
4) The map 
$H^1(G_{S \cup S'  \cup \{s_2\} \backslash \{s'\} },Ad^0(\rhobar))
  \rightarrow \oplus_{v \in S \cup S' \backslash \{s'\}}
  H^1(G_v,Ad^0(\rhobar))/{\cal N}_v$ has
dimensional kernel spanned by $\psi$  
\end{prop}
{\bf Proof:} By the first two remarks at the end of the 
section, 
the Selmer
group for $S \cup S'$ is one dimensional
as is the dual Selmer group for $S \cup S'$.
Let $\psi$ and $\phi$ span the  Selmer
group and dual Selmer group for $S \cup S' \backslash
\{s'\}$ respectively.
Use Lemma \ref{disjointness}
to choose a prime $s_1$ such that
\begin{itemize} 
   \item $1$) is satisfied,
   \item $\psi |_{G_{s_1}}=0$
   \item  $\phi |_{G_{s_1}} \neq 0$.
\end{itemize}
This choice insures that the Selmer group and dual Selmer
group for $S \cup S' \cup \{s_1\} \backslash \{s'\}$
are one dimensional and that (the inflation of) 
$\psi$ spans this Selmer group. 
The map $3$) is surjective by Proposition $1.6$
of [W]. (We take ${\cal L}_{s_1}$
to be all of $H^1(G_{s_1},Ad^0(\rhobar))$
and for $v \neq s_1$ we take
${\cal L}_v$ to be  ${\cal N}_v$).
That the kernel of $3$) is one dimensional also
follows from Proposition $1.6$ of [W].

 Let $\tilde{\phi}$
span the Selmer group for $S \cup S' \cup \{s_1\}
\backslash \{s'\}$.
Now choose $s_2$ to satisfy the three bulleted
items below.
\begin{itemize}
  \item $2$) is satisfied,
  \item  $\psi |_{G_{s_2}}$ is not zero in 
         $H^1(G_{s_2}/I_{s_2}, Ad^0(\rhobar))$,
  \item  $\phi |_{G_{s_2}} \neq 0$,
  \item  $\tilde{\phi} |_{G_{s_2}} \neq 0$.
\end{itemize} 
By construction  $\phi$ and
$\tilde{\phi}$ are independent so Lemma \ref{disjointness} 
allows us to choose $s_2$ as above.
This gives that the sets $S'
\cup \{ s_2\} \backslash \{s'\}$ and $S' 
\cup \{ s_1, s_2\} \backslash \{s'\}$ are both auxiliary.
The Selmer group map for $S \cup S' \cup \{s_2\} \backslash \{s'\}$
is then an isomorphism.
An application of Proposition $1.6$
of [W] gives us $4$) and the proof
is done.

\vspace{3mm}

We have that
\begin{itemize}
\item The sets $S'$, $S' \cup \{s_2\} \backslash \{s'\}$,
and
$S''=S' \cup \{s_1,s_2\} \backslash \{s'\}$ 
are
all auxiliary. 
\item
$\rho_{S \cup S',n-1}^{S'-new}$ is minimal at all
primes of $S$, special at all primes in $S''$,
and unramified outside $S \cup S''$.
Since $R^{S'-new}_{S \cup S'}$ and 
$R^{S''-new }_{S \cup S''}$
are both isomorphic to $W({\bf k})$ we see the deformations
$\rho_{S \cup S'',n-1}^{S''-new}$ 
and $\rho_{S \cup S',n-1}^{S'-new}$ are equal.

\item 
$\rho_{S \cup S',n}^{S'-new}$ is unramified at $s'$
and {\em not} special at $s_1$.
\end{itemize}
Using these facts it is routine to see 
there is an $h \in H^1(G_{S \cup S''}  
,Ad^0(\rhobar))$ such that
$\rho_{S \cup S'',n}^{s^{''}-new}
=(I+p^{n-1}h)
\rho_{S \cup S',n}^{S'-new}$.

\begin{prop}\label{sjiram} $\rho_{S \cup S''}^{S''-new}$
is ramified at both $s_1$ and $s_2$.
\end{prop}
{\bf Proof.} 
We have
$$ \rho_{S \cup S'',n}^{S''-new} 
=(I+p^{n-1}h)
\rho_{S \cup S',n}^{S'-new}.$$
Note that $h |_{G_v} \in {\cal N}_v$ for all
$v \in S \cup S' \backslash \{s'\}$ since 
both sides are the mod $p^n$ reductions
of elements of ${\cal C}_v$ for these $v$.

Suppose the left hand side is unramified
at $s_i$ for $i=1$ or $2$.
Then $h$ inflates from
$H^1(G_{S \cup S'' \backslash \{ s_i\}},
Ad^0(\bar{\rho}))$ and is in the kernel of 
the map $4$) or $3$) of Proposition \ref{h1surj}. But
the kernels of  these  maps are spanned by 
(the inflation of) $\psi$
which is trivial at $s_1$. Then the right hand side
above is not special at $s_1$, so the left hand side is 
also not special, a contradiction. Thus
$\rho_{S \cup S'',n}^{S''-new}$ is
ramified at $s_1$ and $s_2$.

Now replace the auxiliary set $S'$ by the new auxiliary
set $S''=S'
\cup \{s_1,s_2\} \backslash \{s'\}$ and repeat the process.
The fact that there are infinitely many ramified auxiliary sets
follows from an inspection of the proof. 
This finishes the proof of Theorem \ref{main1}.

\section{Selmer groups}

Consider a  ramified auxiliary set $S^{'}=\{s'_1,\dots,s'_n\}$.
Then the corresponding
representation 
$\rho^{S^{'}-new}_{S\cup S^{'}}:
G_{\mathbb{Q}} \rightarrow GL_2(W({\bf k}))$, 
that arises from the isomorphism 
$R_{S \cup S^{'}}^{S^{'}-new} \simeq W({\bf k})$,
is ramified at all the primes in $S^{'}$. Such an $S^{'}$
exists by Theorem \ref{main1}. We have a
surjection 
$R_{S \cup S^{'}} \rightarrow 
R_{S\cup S^{'}}^{S^{'}-new} \simeq W({\bf k})$:
let
$\wp_{S^{'}}$ be the kernel. For each $i$ let $m_i$ be
the largest integer such that $\rho_{S\cup S^{'},m_i}^{S'-new}|_{G_{s'_i}}$ 
is unramified. By choice we have $m_i < \infty$.

\begin{theorem}\label{finiteness1}
  If ${S^{'}}$ is a ramified auxiliary set,
  the tangent space $\wp_{S^{'}}/\wp_{S^{'}}^2$ is a finite abelian group.
\end{theorem}
{\bf Proof.} We see by standard arguments as in
Proposition 1.2 of [W], or [Ma1], 
that we can identify the dual of $\wp_{S^{'}}/\wp_{S^{'}}^2$ 
as a Selmer group, i.e., as
the kernel of the map 
\begin{eqnarray*}
  H^1(G_{S \cup {S^{'}}},Ad^0(\rho_{S \cup {S^{'}}}^{{S^{'}}-new})\otimes 
  {\mathbb{Q}}_p/{\mathbb{Z}}_p) \rightarrow \Pi_{v \in S }
  {{H^1(G_v, Ad^0(\rho_{S \cup {S^{'}}}^{{S^{'}}-new})\otimes 
  {\mathbb{Q}}_p/{\mathbb{Z}}_p)} \over {L_v}}.
\end{eqnarray*}
For $s'_i \in {S^{'}}$, note the representation space
of $\rho_{S \cup {S^{'}}}^{{S^{'}}-new}|_{G_{s'_i}}$
has a one dimensional stable subspace. Let $Z \subset Ad^0
(\rho_{S \cup {S^{'}}}^{{S^{'}}-new})$ be the (upper triangular)
one dimensional subspace of nilpotent
elements that are preserved under conjugation by
$\rho_{S\cup {S^{'}}}^{{S^{'}}-new}|_{G_{s'_i}}$.

We define $L_v$ as 
follows:
\begin{itemize}
    \item For $v \in S$, $v \neq p$, we define $L_v$
          to be the kernel of the map $$H^1(G_v,
          Ad^0(\rho_{S \cup {S^{'}}}^{{S^{'}}-new})\otimes
          \mathbb{Q}_p/\mathbb{Z}_p) \rightarrow
          H^1(I_v,
          Ad^0(\rho_{S \cup {S^{'}}}^{{S^{'}}-new})\otimes
          \mathbb{Q}_p/\mathbb{Z}_p).$$
    
    \item For $s' \in S'$ we define $L_{s'}$ to be the kernel of 
          the map 
$$H^1(G_{s'},Ad^0(\rho_{S \cup {S^{'}}}^{{S^{'}}-new})
\otimes {\mathbb{Q}}_p/{\mathbb{Z}}_p)
\rightarrow H^1(G_{s'},(Ad^0(\rho_{S \cup {S^{'}}}^{{S^{'}}-new})/Z)
\otimes {\mathbb{Q}}_p/{\mathbb{Z}}_p).$$
    \item If $\rhobar |_{G_p}$ is odd, reducible, and flat, $L_p$ is
          defined as in the flat case.
    \item If $\rhobar |_{G_p}$ is odd and 
          not reducible then we define $L_p$ as in the flat
          case using Fontaine-Lafaille modules of appropriate
          filtration length.
    \item In all other (necessarily reducible)
          odd cases define $L_p$ as in the ordinary 
          situation.
    \item If  $\rhobar$ is even then $L_p$ is defined to be all
          of $H^1(G_p,Ad^0(\rho_{S \cup {S^{'}}}^{{S^{'}}-new})\otimes
          {\mathbb{Q}}_p/{\mathbb{Z}}_p)$.
\end{itemize}
Recall we identify 
the dual of $\wp_{{S^{'}}}/\wp_{{S^{'}}}^2$
with the kernel of the map
\begin{eqnarray*}
  H^1(G_{S \cup {S^{'}}},Ad^0(\rho_{S \cup {S^{'}}}^{{S^{'}}-new})\otimes 
  {\mathbb{Q}}_p/{\mathbb{Z}}_p) \rightarrow \Pi_{v \in S }
  {{H^1(G_v, Ad^0(\rho_{S \cup {S^{'}}}^{{S^{'}}-new})\otimes 
  {\mathbb{Q}}_p/{\mathbb{Z}}_p)} \over {L_v}}.
\end{eqnarray*}
Using the fact that $R_{S \cup {S^{'}}}^{{S^{'}}-new}$ 
is isomorphic to $W({\bf k})$, we deduce that
the kernel of the map
\begin{eqnarray*}
  H^1(G_{S \cup {S^{'}}},Ad^0(\rho_{S \cup {S^{'}}}^{{S^{'}}-new})\otimes 
  {\mathbb{Q}}_p/{\mathbb{Z}}_p) \rightarrow \Pi_{v \in S \cup {S^{'}}}
  {{H^1(G_v, Ad^0(\rho_{S \cup Q}^{Q-new})\otimes 
  {\mathbb{Q}}_p/{\mathbb{Z}}_p)} \over {L_v}}
\end{eqnarray*}
is trivial. 
Thus the dual of $\wp_{S^{'}}/\wp_{S^{'}}^2$ injects into 
$$ \Pi_{s' \in  {S^{'}}} 
{{H^1(G_{s'}, 
Ad^0(\rho_{S \cup {S^{'}}}^{{S^{'}}-new})\otimes
  {\mathbb{Q}}_p/{\mathbb{Z}}_p)}   \over {L_{s'}}}.$$
The following lemma completes the proof of Theorem 
\ref{finiteness1}.

\begin{lemma}\label{localselmer}
 Let $s' \in {S^{'}}$ and let $m$ be the largest integer such that
  $\rho_{S \cup {S^{'}},m}^{{S^{'}}-new} |_{G_{s'}}$ is unramified. 
\newline\noindent
1) $(Ad^0(\rho_{S \cup {S^{'}}}^{{S^{'}}-new})\otimes
  {\mathbb{Q}}_p/{\mathbb{Z}}_p)^{I_{s'}} \supseteq
  Ad^0(\rho_{S \cup {S^{'}}}^{{S^{'}}-new})\otimes
   p^{-m}\mathbb{Z}/\mathbb{Z}$ and
   $${ {   (Ad^0(\rho_{S \cup S'}^{S'-new})\otimes
  \mathbb{Q}_p/\mathbb{Z}_p)^{I_{s'}} } \over
  { Ad^0(\rho_{S \cup S'}^{S'-new})\otimes
   p^{-m}\mathbb{Z}/\mathbb{Z}   }} \simeq Z \otimes
     \mathbb{Q}_p/\mathbb{Z}_p.$$
\newline\noindent
2) The inflation map  
  $$H^1(G_{s'}/I_{s'}, 
  (Ad^0(\rho_{S \cup {S^{'}}}^{{S^{'}}-new})\otimes
  {\mathbb{Q}}_p/{\mathbb{Z}}_p)^{I_{s'}}) \rightarrow
  H^1(G_{s'}, (Ad^0(\rho_{S \cup {S^{'}}}^{{S^{'}}-new})\otimes
  {\mathbb{Q}}_p/{\mathbb{Z}}_p))$$ is an isomorphism.
\newline\noindent
3)  $ H^1(G_{s'}, Ad^0(\rho_{S \cup {S^{'}}}^{{S^{'}}-new})\otimes
  {\mathbb{Q}}_p/{\mathbb{Z}}_p) \simeq W({\bf k})/p^{m}$
\newline\noindent
4) For all $s' \in {S^{'}}$ we have that $L_{s'}$ is trivial
\end{lemma}
{\bf Proof.} 
First observe since $s'$ is {\em not} $1$ mod $p$ that 
$H^r(G_{s'}/I_{s'},Z \otimes
  {\mathbb{Q}}_p/{\mathbb{Z}}_p)$ and
$H^r(G_{s'}/I_{s'},Z \otimes p^{-m}\mathbb{Z}/\mathbb{Z})$ are
trivial for $r=0,1$.
\newline\noindent
1) This is a routine computation that uses ramification
first occurs mod $p^{m+1}$.
\newline\noindent
2) Using inflation-restriction,
it suffices to prove $H^1(I_{s'}, 
Ad^0(\rho_{S \cup {S^{'}}}^{{S^{'}}-new})\otimes
  {\mathbb{Q}}_p/{\mathbb{Z}}_p)$ is trivial.
This last group is isomorphic to
$( Ad^0(\rho_{S \cup {S^{'}}}^{{S^{'}}-new})\otimes
  {\mathbb{Q}}_p/{\mathbb{Z}}_p)_{I_{s'}}$ which in turn
is isomorphic to $W({\bf k})(-1) \otimes
  {\mathbb{Q}}_p/{\mathbb{Z}}_p$ which is trivial as $s'$ is
{\em not} $1$ mod $p$.   
\newline\noindent
3)
We have the exact sequence
$$0 \rightarrow Ad^0(\rho_{S \cup S^{'}}^{{S^{'}}-new})\otimes
p^{-m}{\mathbb{Z}}/{\mathbb{Z}} \rightarrow
(Ad^0(\rho_{S \cup {S^{'}}}^{{S^{'}}-new})
\otimes {\mathbb{Q}}_p/{\mathbb{Z}}_p)^{I_{s'}}
\rightarrow Z \otimes {\mathbb{Q}}_p/{\mathbb{Z}}_p \rightarrow 0$$
of $G_{s'}/I_{s'}$ modules.
Since $s'$ is not $\pm 1$ mod $p$ we see
$H^i(G_{s'}/I_{s'}, Z \otimes \mathbb{Q}_p/\mathbb{Z}_p)=0$
for $i=0,1$.
Taking $G_{s'}/I_{s'}$ cohomology
 we conclude
$$H^1(G_{s'}/I_{s'}, (Ad^0(\rho_{S \cup {S^{'}}}^{{S^{'}}-new})
\otimes {p^{-m}\mathbb{Z}}/{\mathbb{Z}}))
\simeq H^1(G_{s'}/I_{s'},
(Ad^0(\rho_{S \cup {S^{'}}}^{{S^{'}}-new})
\otimes {\mathbb{Q}}_p/{\mathbb{Z}}_p)^{I_{s'}}).$$

Since
$$Ad^0(\rho_{S \cup {S^{'}}}^{{S^{'}}-new})
\otimes p^{-m}{\mathbb{Z}}/{\mathbb{Z}}
\simeq
W({\bf k})(-1)/p^m \oplus W({\bf k})/p^m \oplus W({\bf k})(1)/p^m$$ as
$G_{s'}/I_{s'}$ modules, the conclusion
follows from $2$) and a routine calculation.
\newline \noindent
4)
We take $G_{s'}$ cohomology of the sequence
$$0 \rightarrow Z \otimes {\mathbb{Q}}_p/{\mathbb{Z}}_p
\rightarrow Ad^0(\rho_{S \cup {S^{'}}}^{{S^{'}}-new})
\otimes {\mathbb{Q}}_p/{\mathbb{Z}}_p
\rightarrow Ad^0(\rho_{S \cup {S^{'}}}^{{S^{'}}-new})/Z
\otimes {\mathbb{Q}}_p/{\mathbb{Z}}_p \rightarrow 0.$$
The computation of $G_{s'}$ invariants is routine
and we get
$$0 \rightarrow 0 \rightarrow W({\bf k})/p^m \rightarrow
W({\bf k}) \otimes {\mathbb{Q}}_p/{\mathbb{Z}}_p
\rightarrow H^1(G_{s'}, Z \otimes {\mathbb{Q}}_p/{\mathbb{Z}}_p)$$
$$\rightarrow  H^1(G_{s'},Ad^0(\rho_{S \cup {S^{'}}}^{{S^{'}}-new})
\otimes {\mathbb{Q}}_p/{\mathbb{Z}}_p)
\rightarrow  H^1(G_{s'},Ad^0(\rho_{S \cup {S^{'}}}^{{S^{'}}-new})/Z
\otimes {\mathbb{Q}}_p/{\mathbb{Z}}_p) \rightarrow ...$$
We will prove the map $W({\bf k}) \otimes {\mathbb{Q}}_p/{\mathbb{Z}}_p
\rightarrow H^1(G_{s'}, Z \otimes {\mathbb{Q}}_p/{\mathbb{Z}}_p)$
is surjective. This will give injectivity of
$$H^1(G_{s'},Ad^0(\rho_{S \cup {S^{'}}}^{{S^{'}}-new})
\otimes {\mathbb{Q}}_p/{\mathbb{Z}}_p)
\rightarrow  H^1(G_{s'},Ad^0(\rho_{S \cup {S^{'}}}^{{S^{'}}-new})/Z
\otimes {\mathbb{Q}}_p/{\mathbb{Z}}_p)$$ and we will be done.

Recall $(Z \otimes {\mathbb{Q}}_p/{\mathbb{Z}}_p)^{I_{s'}}
=Z \otimes {\mathbb{Q}}_p/{\mathbb{Z}}_p$ and
consider the inflation-restriction sequence
$$0 \rightarrow 
H^1(G_{s'}/I_{s'},(Z \otimes {\mathbb{Q}}_p/{\mathbb{Z}}_p)^{I_{s'}})
\rightarrow H^1(G_{s'},Z \otimes {\mathbb{Q}}_p/{\mathbb{Z}}_p)$$
$$\rightarrow H^1(I_{s'},Z \otimes {\mathbb{Q}}_p/{\mathbb{Z}}_p)^{G_{s'}
/I_{s'}} \rightarrow 0.$$
The last (fifth) term is $0$ as 
it is actually  an $H^2$ of the group
$G_{s'}/I_{s'}$
which has cohomological dimension $1$.
We know
$H^1(G_{s'}/I_{s'},(Z \otimes {\mathbb{Q}}_p/{\mathbb{Z}}_p)^{I_{s'}})$
is trivial so
$$H^1(G_{s'},Z \otimes {\mathbb{Q}}_p/{\mathbb{Z}}_p)
\rightarrow H^1(I_{s'},Z \otimes {\mathbb{Q}}_p/{\mathbb{Z}}_p)^{G_{s'}
/I_{s'}}$$
is an isomorphism.
Since $I_{s'}$ acts trivially on $Z$ we see
$ H^1(I_{s'},Z \otimes {\mathbb{Q}}_p/{\mathbb{Z}}_p)^{G_{s'}
/I_{s'}}$ is just the $G_{s'}/I_{s'}$ equivariant
elements of $Hom(I_{s'},Z \otimes {\mathbb{Q}}_p/{\mathbb{Z}}_p)$.
By our choice of $s'$ the actions of  $G_{s'}/I_{s'}$ 
on (the abelianisation of)
$I_{s'}$ and on $Z \otimes {\mathbb{Q}}_p/{\mathbb{Z}}_p)$ are
compatible so we need only find
$Hom(I_{s'},Z \otimes {\mathbb{Q}}_p/{\mathbb{Z}}_p) $
which is just 
$Hom( \mathbb{Z}_p, W({\bf k}) \otimes {\mathbb{Q}}_p/{\mathbb{Z}}_p)$
which in turn is $W({\bf k}) \otimes {\mathbb{Q}}_p/{\mathbb{Z}}_p$.

We have proved $H^1(I_{s'},Z \otimes {\mathbb{Q}}_p/{\mathbb{Z}}_p)^{G_{s'}
/I_{s'}}$ is isomorphic to
$W({\bf k}) \otimes {\mathbb{Q}}_p/{\mathbb{Z}}_p$.
Thus the map we needed to prove surjective is just the multiplication
by $p^m$ map on $W({\bf k}) \otimes {\mathbb{Q}}_p/{\mathbb{Z}}_p$ which
is clearly onto.
The proposition and Theorem \ref{finiteness1} are now proved.

\vspace{3mm}

\noindent{\bf Remarks:}
\newline\noindent
$1$) Note that the proof of Therorem
\ref{finiteness1} crucially uses, and is more or less
a formal consequence of the fact, that
$\rho_{S \cup {S^{'}}}^{{S^{'}}-new}|_{G_{s'_i}}$ is 
ramified for all $s'_i \in {S^{'}}$.
\newline\noindent
$2$) It also follows from $3$) and the proof that if we denote by
$\wp_{\alpha}$ the kernel of the map $R_{S \cup {S^{'}}}^{\alpha-new} 
\rightarrow  R_{S \cup {S^{'}}}^{{S^{'}}-new}$, for any 
$\alpha \subset {S^{'}}$,
then $\wp_{\alpha}/{\wp_{\alpha}}^2$ is isomorphic to a subgroup of
$\Pi_{s'_i \notin \alpha}W({\bf k})/p^{m_i}$.

\vspace{3mm}

\begin{cor}\label{minimal1}
  The prime ideal $\wp_{S^{'}}$ 
is minimal and thus ${\rm Spec}(R_{S \cup {S^{'}}})$
  has an isolated component that is isomorphic to $W({\bf k})$.
\end{cor}

\noindent{\bf Proof.}
Consider the localisation $R_{\wp_{S^{'}}}$ of 
$R_{S \cup {S^{'}}}$ at $\wp_{S^{'}}$.
Then as $R_{S \cup {S^{'}}}/\wp_{S^{'}}$ is 
of characteristic 0, we see that
$p$ is invertible in $R_{\wp_{S^{'}}}$, and 
thus that $\wp_{S^{'}}/\wp_{S^{'}}^2$ is zero
as it follows from Theorem \ref{finiteness1} 
that it is a $p$-torsion group.
But then as $\wp_{S^{'}}$ is a proper ideal, and $R_{\wp_{S^{'}}}$ 
is Noetherian, we conclude 
that $\wp_{S^{'}}$ is 0 in $R_{\wp_{S^{'}}}$. Hence 
$\wp_{S^{'}}$ is a minimal
prime ideal of $R_{S \cup {S^{'}}}$, i.e., 
${\rm Spec}(R_{S \cup {S^{'}}}^{{S^{'}}-new})$
is an isolated component of ${\rm Spec}(R_{S \cup {S^{'}}})$. 

\begin{cor}
  If $R_{S \cup {S^{'}}}$ is a complete intersection ring, 
  then it is finite flat
  over $W({\bf k})$.
\end{cor}

\noindent{\bf Proof.} We can arrange for an epimorphism
$\phi:W({\bf k})[[T_1,\dots,T_d]] \rightarrow R_{S \cup {S^{'}}}$
with kernel $I$ contained in the ideal
$J=(T_1,\dots,T_d)$, and further the composite map
$W({\bf k})[[T_1,\dots,T_d]] \rightarrow R_{S \cup {S^{'}}}
\rightarrow R_{S \cup {S^{'}}}^{{S^{'}}-new}$ 
given by $T_i \rightarrow 0$. Let $I$ be generated
by $(r_1,\dots,r_{e})$ for some integer $e \geq 0$. 
As $\phi(J)/\phi(J)^2$ is finite by Theorem \ref{finiteness1}
the images of $r_i$ in $J/J^2$ span
a full sublattice of $J/J^2 \simeq W({\bf k})^d$. 
Thus $e$ is at least $d$ and the corollary follows.

\vspace{3mm}

 Let   $\gamma_{s'}$ be a generator
of the unique ${\mathbb{Z}}_p$ quotient of $I_{s'}$.
Consider the universal representation 
$\rho_{S \cup {S^{'}}}^{univ}:G_{\mathbb{Q}} 
\rightarrow GL_2(R_{S \cup {S^{'}}})$. 
Because of the nature of the primes
$s' \in {S^{'}}$ (in particular none of them is
congruent to $\pm 1$ mod $p$) 
and the structure of tame inertia groups recalled above
we see that 
$\rho_{S \cup {S^{'}}} ^{univ}({\gamma_{s'}})$ is well-defined 
(on choosing a prime of ${\bf{\bar Q}}$ above $s'$) and of the form 
$$\left(\matrix{1&x_{s'}\cr
                0&1\cr}\right).$$ 
Consider the principal ideal $(\Pi_{s' \in {S^{'}}}x_{s'})$
and choose a generator 
$x_{S^{'}}$ for it. Note that $x_{S^{'}}$ is not nilpotent
as ${S^{'}}$ is a ramified auxiliary set.

\begin{cor}
  $R_{S \cup {S^{'}}}[{1 \over x_{S^{'}}}] \simeq Frac(W({\bf k}))$.
\end{cor}

\noindent{\bf Proof.} Consider the representation
$G_{\mathbb{Q}} \rightarrow GL_2(R_{S \cup {S^{'}}}[{1 \over x_{S^{'}}}])$,
the composition of $\rho_{S \cup {S^{'}}}^{univ}$ with
the map $R_{S \cup {S^{'}}} 
\rightarrow R_{S \cup {S^{'}}}[{1 \over x_{S^{'}}}]$. In this
representation we deduce from the structure of tame inertia,
and the fact that $x_{S^{'}}$ is not a zero-divisor in
$R_{S \cup {S^{'}}}[{1 \over x_{S^{'}}}]$, 
that for $s' \in {S^{'}}$ the ratio of the eigenvalues
of the image of a lift of ${\rm Frob}_{s'}$ is $s'$. 
Now upon using the facts that $R_{S \cup {S^{'}}}^{{S^{'}}-new}=W({\bf k})$ 
and $0 \neq (x_{S^{'}}) \subset (p)$ in $R^{{S^{'}}-new}_{S \cup S^{'}}$,
the corollary follows.

\vspace{3mm}

We end this section by indicating an application of
Theorem \ref{finiteness1} to a conjecture of Fontaine-Mazur. 
Conjecture 2c of [FM] in our situation translates as
saying that there are only finitely many morphisms 
$R_{S \cup {S^{'}}} \rightarrow W({\bf k})$ (taking $K=W({\bf k})$ in
the notation of loc. cit.). By a use of the finiteness theorem
of Hermite-Minkowski this is equivalent to saying that 
a morphism $\pi:R_{S \cup {S^{'}}} \rightarrow W({\bf k})$ cannot be the
nontrivial limit of morphisms $\pi_i:R_{S \cup {S^{'}}} \rightarrow 
W({\bf k})$, i.e., if for any $n$, $\pi_i \equiv \pi \pmod{p^n}$
for $i$ large enough, then $\pi_i$ is eventually constant. 
In this situation we provide an application of Theorem \ref{finiteness1}:

\begin{prop}\label{fm}
  If $\pi:R_{S \cup {S^{'}}} 
  \rightarrow W({\bf k})$ and ${\rm ker}(\pi)/{\rm
  ker}(\pi)^2$ is a finite abelian group, then any sequence of morphisms
  $\pi_i:R_{S \cup {S^{'}}} \rightarrow 
  W({\bf k})$ that tends to $\pi$ is eventually constant. Thus for
  a ramified auxiliary set ${S^{'}}$, 
  the morphism $R_{S \cup {S^{'}}} \rightarrow
  R_{S \cup {S^{'}}}^{{S^{'}}-new} \simeq W({\bf k})$ 
  cannot be the non trivial limit of morphisms
  $\pi_i:R_{S \cup {S^{'}}} \rightarrow W({\bf k})$.
\end{prop}

\noindent{\bf Proof:} The second part follows from Theorem
\ref{finiteness1} and the first. We prove the first part. 
Let $p^n$ be the exponent of ${\rm ker}(\pi)/{\rm ker}(\pi)^2$ and we
claim that if $\pi_i \equiv \pi$ mod $p^{n+1}$, then $\pi_i= \pi$.
If not there is a $n' \geq n+1$ such that $\pi_i \equiv \pi$ mod
$p^{n'}$ but $\pi_i$ is not congruent to $\pi$  mod $p^{n'+1}$. But then 
$\pi_i - \pi:\wp/\wp^2 \rightarrow p^{n'}W({\bf k})/p^{2n'}W({\bf k})$
is surjective. This contradicts the fact that $p^n$ is the 
exponent of ${\rm ker}(\pi)/{\rm ker}(\pi)^2$.

\section{Odds and ends}

\subsection{Parity}

We prove a proposition that is needed in [K].

\begin{prop}\label{parity}
  We can choose an auxiliary set ${\tilde{S}}$ 
  so that its cardinality is of any given parity.
\end{prop}

\noindent{\bf Proof.} We sketch the proof. 
Consider any auxiliary set ${S^{'}}$. 
We may assume that it is not of the 
desired parity. We add a prime $w_1$ not congruent to $\pm 1$ mod $p$
that is special for $\rhobar$. 
If ${S^{'}} \cup \{w_1\}$ is an auxiliary set, take
${\tilde{S}} = {S^{'}} \cup \{w_1\}$ and we are
done. Otherwise the Selmer group 
$H^1_{ {\cal N}_v}(G_{S \cup {S^{'}} \cup \{w_1\}},Ad^0(\rhobar))$ 
is  one dimensional as is the 
dual Selmer group 
$H^1_{ {\cal N}^{\perp}_v}
(G_{S \cup {S^{'}} \cup \{w_1\}},Ad^0(\rhobar)^*)$.
Let $\psi$ span the Selmer group and
$\phi$ span the dual Selmer group.
Simply choose a prime $w_2$ such that $\psi$ is trivial
at $G_{w_2}$ and $\phi$ does not map to zero in 
$H^1(G_{w_2},Ad^0(\rhobar)^*)$.
We then see that the dual Selmer group
$H^1_{ {\cal N}^{\perp}_v}(G_{S \cup {S^{'}} \cup \{w_1,w_2\}},
Ad^0(\rhobar)^*)$
is one dimensional. Now annihilate this dual Selmer group by choosing
an appropriate $w_3$ and take $\tilde{S}$
to be the  auxiliary set ${S^{'}} \cup \{w_1,w_2,w_3\}$
and we are done.

\subsection{Finer structure of deformation rings}

We  make the following conjecture for any ramified
auxiliary set $S^{'}=\{s'_1,\dots ,s'_n\}$ 
that is guaranteed to exist by Theorem \ref{main1}.

\begin{conjecture}\label{flatness}
  The ring $R_{S \cup {S^{'}}}$ is a finite flat complete
  intersection over $W({\bf k})$.
\end{conjecture}

\vspace{3mm}

\noindent{\bf Remark:}
\newline\noindent
 The recent work of Taylor, cf. [T1], will imply this for odd $\rhobar$
  in many cases (using G.~B\"ockle's arguments in the 
  appendix to [K]). There is some
  evidence in the even case in [R4] and [B].

\begin{conjecture}\label{selmerstructure}
  In the situation of Theorem \ref{finiteness1}, the abelian group
  $\wp_{S^{'}}/\wp_{S^{'}}^2$ that was proved to be a subgroup of
  $\oplus_1^nW({\bf k})/(p^{m_i})$ in 
  Theorem \ref{finiteness1}, 
  is in fact isomorphic to $\oplus_1^nW({\bf k})/(p^{m_i})$
\end{conjecture}

Let $\sigma_{{s'_i}}$ be a lift of Frobenius in the Galois
group of the maximal tame extension of ${\mathbb{Q}}_{{s'_i}}$ as before.
Denote by $\alpha_{{s'_i}}$ the ratio of the eigenvalues
of $\rho_{S^{'}}^{univ}(\sigma_{{s'_i}})$. The $\alpha_{{s'_i}}$ are well-defined and
independent of choice of the lift $\sigma_{{s'_i}}$.
Using the results of [Bo], and because of our assumption that ${s'_i}$ is
not $\pm 1$ mod $p$ we have:

\begin{prop}\label{boston}
  The kernel of the surjective map $R_{S \cup S^{'}} 
  \rightarrow R_{S \cup
  S^{'}}^{S^{'}-new}$ is the ideal generated by 
  $(\alpha_{s'_1}-s'_1,\dots,\alpha_{s'_n}-s'_n)$.
\end{prop}

\subsection{Infinite ramification}

In [R1] it was proved that there exist odd $2$ dimensional
continuous, 
irreducible $p$-adic Galois representations 
ramified at infinitely many primes.
Theorem $2$ of that paper proved, assuming the GRH, the
existence of odd $2$ dimensional 
$p$-adic Galois representations ramified at infinitely 
many primes that were crystalline at $p$.
Here we remove the GRH hypothesis  and prove there exist even
$2$ dimensional $p$-adic Galois representations
ramified at infinitely many primes.

\begin{theorem}\label{infram}  Let 
$\rhobar:G_{\mathbb{Q}}\rightarrow GL_2({\bf k})$ 
satisfy the hypotheses of the introduction.
\newline
a) If ${\bar {\rho}}$ is even 
then there exists a
deformation 
${\rho}:G_{\mathbb{Q}}  \rightarrow GL_2(W({\bf k}))$ of  $\rhobar$
ramified at infinitely many primes.
\newline b) If ${\bar {\rho}}$ is odd 
then there exists a a
deformation 
${\rho}:G_{\mathbb{Q}} \rightarrow GL_2(W({\bf k}))$ of  $\rhobar$
ramified at infinitely many primes that is potentially semistable
at $p$.
\end{theorem}

\noindent{\bf Remark:}
\newline\noindent 
For $E_{/{\mathbb{Q}}}$ an elliptic curve without
complex multiplication it is a well known result of Serre
that for almost all $p$, the $p$-torsion of $E$ gives rise
to a ${\bar {\rho}}$ to which part {\em b)} of the theorem applies.

\noindent{\bf Proof.}
Note for ${\bar {\rho}}$ in part b) that ${\cal C}_p$
consists of either 
\begin{itemize}
    \item deformations of
          ${\bar {\rho}}|_{G_p}$ of fixed positive rational
          integral weight that are after a fixed finite twist ordinary
    \item or deformations that after a fixed finite twist come from
          irreducible Fontaine-Lafaille modules of fixed
          filtration length.
\end{itemize}
We will construct our global representation as a $p$-adic
limit of global representations whose restrictions
to $G_p$ all lie in ${\cal C}_p$. Since each of the two
classes above is closed under limits we get that
our limit representation is either potentially ordinary 
or potentially crystalline, hence potentially semistable.

Choose an auxiliary set $S'_0$ for $\rhobar$.
We have 
deformation $\rho^{S'_0 - new}_{S \cup S'_0}$ 
of $\rhobar$ corresponding to the isomorphism
$R^{S'_0 - new}_{S \cup S'_0} \simeq W({\bf k})$.
We proceed to inductively construct a sequence of characteristic
zero representations ramified at more and more primes whose $p$-adic
limit is the desired $\rho$.

Let $n=1$.
As in Lemma \ref{disjointness} choose a prime $q_1$ 
unramified in $\rho^{S'_0 - new}_{S \cup S'_0}$  and   
$\rho^{S'_0-new}_{S \cup S'_0,n}$
is special at $q_1$ (in particular
$q_1 \not \equiv \pm 1$ mod $p$) but 
$\rho^{S'_0-new}_{S \cup S'_0,n+1}$ is not special at $q_1$.
If $S'_0 \cup \{q_1\}$ is auxiliary
then there is the 
unique deformation $\rho^{S'_0 \cup \{q_1\}-new}_{S \cup S'_0 \cup \{q_1\}}$ 
of ${\bar {\rho}}$ 
corresponding to the isomorphism
$R^{S'_0 \cup \{q_1\} - new}_{S \cup S'_0 \cup \{q_1\}} \simeq W({\bf k})$.
This unique deformation cannot be unramified at $q_1$.
For if it were, the eigenvalues of
$\rho^{S'_0 \cup \{q_1\}-new}_{S \cup S'_0 \cup \{q_1\}}(Frob_{q_1})$ 
would have ratio $q_1$. Being unramified at $q_1$
this characteristic
zero deformation of ${\bar {\rho}}$ would necessarily be 
$\rho^{S'_0-new}_{S \cup S'_0}$. But
the eigenvalues of $\rho^{S'_0-new}_{S \cup S'_0}(Frob_{q_1})$ 
do {\em not}
have ratio $q_1$ by choice.

Suppose then that 
$S'_0 \cup \{q_1\}$ is not an auxiliary set.
Following Lemma \ref{disjointness}
we can find a  prime $q_2$
such that:
\begin{itemize}
\item $\rho^{S'_0 \cup \{q_1\} - new}_{S \cup S'_0 \cup \{q_1\},n}$ 
is special at $q_2$ but
$\rho^{S'_0 \cup \{q_1\} - new}_{S \cup S'_0 \cup \{q_1\},n+1}$ 
is not special at $q_2$. 
\item  $ S''_0=S'_0 \cup \{q_1,q_2\}$ is auxiliary.
\end{itemize}
Then 
$\rho^{S''_0  - new}_{S \cup S''_0 }$
is the  unique deformation of ${\bar {\rho}}$ 
corresponding to the isomorphism
$R^{S''_0  - new}_{S \cup S''_0 } 
\simeq W({\bf k})$.
It is ramified at at least one of $q_1,q_2$,
for otherwise we would have 
$\rho^{S''_0 - new}_{S \cup S''_0} 
=\rho^{S'_0-new}_{S \cup S'_0}$
and $\rho^{S'_0-new}_{S \cup S'_0}(Frob_{q_1})$ and  
$\rho^{S'_0-new}_{S \cup S'_0}(Frob_{q_2})$ would have eigenvalues
with ratio $q_1$ and $q_2$ respectively, a contradiction to
how these primes were chosen.

If we are in the first case put $S'_1 = S'_0 \cup \{q_1\}$.
If we are in the second case put $S'_1 = S''_0=S'_0 \cup \{q_1,q_2\}$.
In either case $S'_1$ is auxiliary, 
$\rho^{S'_0-new}_{S \cup S'_0,n} = \rho^{S'_1-new}_{S \cup S'_1,n}$
and $\rho^{S'_1-new}_{S \cup S'_1}$ is ramified at 
at least one prime of $S'_1 \backslash S'_0$.
Continuing the induction
we get a sequence of deformations $\rho^{S'_k-new}_{S \cup S'_k}$
such that  $\rho^{S'_k-new}_{S \cup S'_k}$ is ramified at 
at least one prime of $S'_k \backslash S'_{k-1}$
{\em and} this sequence converges $p$-adically. The limit 
deformation is ramified at infinitely many primes.

\vspace{3mm}

\noindent{\bf Remarks:}
\newline\noindent
$1$) If $\rhobar$ is odd and has Serre weight 
      between $2$ and $p$, one sees from the proof
      that we can arrange for the $\rho$
      of the theorem to also be crystalline.
\newline\noindent
$2$) Fontaine and Mazur have conjectured in [FM] that all irreducible
    finite dimensional  
    $p$-adic representations of the absolute Galois group of a number
    field ${\bf F}$ that are potentially semistable at all primes above
    $p$ and have finite ramification
    arise in the \'{e}tale cohomology of some variety over
    ${\bf F}$. Part b) of Theorem \ref{infram} shows that hypotheses
    of finite ramification and potential semistability are independent. 
\newline\noindent
$3$) Assuming ${\bar {\rho}}$ is modular, Theorem 1 of [K] 
   can be used to give another proof of part b) of 
   Theorem \ref{infram}.
\newline\noindent
$4$)  Mestre has shown in [Me] that $SL_2({\bf F_7})$
       is a regular extension of ${\mathbb{Q}}(T)$. 
       Specialization gives us an infinite
       supply of $SL_2({\bf F_7})$ extensions of ${\mathbb{Q}}$.
       It is an exercise in local class field theory to see that
       none of these representations, when restricted to
       $G_7$, are twists of 
       $\left(\begin{array}{cc} \chi &0 \\0 & 1\end{array}\right)$
       or
       the indecomposable
       representation $\left(\begin{array}{cc} \chi^{p-2} & *\\
       0& 1\end{array}\right)$. 
       Thus part a) of Theorem 
       \ref{infram} applies in all these situations.

\section{References}

\noindent [B] B\"ockle, G., {\it A 
local-to-global principle for deformations of Galois
representations}, J. Reine Angew. Math. 509 (1999), 199--236.

\vspace{3mm}

\noindent [Bo] Boston, N., {\it Families of Galois 
representations---increasing the
ramification}, Duke Math. J. 66 (1992), 357--367.

\vspace{3mm}

\noindent [Ca] Carayol, H., 
{\it Formes modulaires et repr\'esentations 
galoisiennes avec valeurs
dans un anneau local complet}, in {\it $p$-adic monodromy 
and the Birch and Swinnerton-Dyer conjecture}, 213--237, 
Contemp. Math., 165, AMS, 1994. 

\vspace{3mm}

\noindent [DDT] Darmon, H., Diamond, F., Taylor, R. 
                {\it Fermat's last theorem} in 
                Elliptic curves, modular forms \& Fermat's
                last theorem (Hong Kong) 1993, p. 2-140. 

\vspace{3mm}
\noindent [Di] Dickson, L. E., {\em Linear Groups}, B. G. Teubner 1901.
\vspace{3mm}

\noindent [FLT] {\it Modular forms and Fermat's last theorem}, 
edited by Cornell, G.,  Silverman, J.H.,  and Stevens, G.
Springer-Verlag, New York, 1997.

\vspace{3mm}

\noindent [FM] Fontaine, J.-M.,  Mazur, B., {\it Geometric Galois 
representations}, Elliptic curves, modular forms, and Fermat's last 
theorem (Hong Kong, 1993), 41--78, Ser. Number Theory, I, 
Internat. Press, Cambridge, MA, 1995. 

\vspace{3mm}

\noindent [Ma1] Mazur, B., {\it An introduction to the deformation
theory of Galois representations}, in [FLT].

\vspace{3mm}

\noindent [Ma2] Mazur, B., {\it Deforming Galois representations}, 
in Galois Groups over $\mathbb{Q}$, eds. Y. Ihara, K. Ribet, J. -P. Serre.

\vspace{3mm}

\noindent[Me]  Mestre, J.-F., {\it Construction d'extensions
                  r\'{e}guli\`{e}res
                  de ${\mathbb{Q}}(T)$ \`{a} groupe de Galois $SL_2({\bf F_7})$ et
                  $\widetilde{M_{12}}$},
                  C.R. Acad. Sci. Paris 319 (1994),  781--782.
\vspace{3mm}

\noindent [Mi]  Milne, J., {\it Arithmetic Duality Theorems},
                  Academic Press, Inc. 1986.
\vspace{3mm}

\noindent [K] Khare, C., {\it On isomorphisms between deformation rings and
                         Hecke rings}, preprint.

\vspace{3mm}

\noindent [R1] Ramakrishna, R., {\it Infinitely ramified representations},
          Annals of Mathematics 151 no. 2 (2000), 793--815.

\vspace{3mm}

\noindent [R2] Ramakrishna, R., {\it Deforming Galois Representations and the
          Conjectures of Serre and Fontaine-Mazur},
          Annals of Mathematics 156 no. 1 (2002), 115--154.

\vspace{3mm}

\noindent [R3] Ramakrishna, R., {\it Lifting Galois Representations},
          Inventiones Mathematicae 138 (1999) 537--562.

\vspace{3mm} 

\noindent [R4] Ramakrishna, R., {\it Deforming an even representation},
Inventiones Mathematicae 132 (1998), 563--580

\vspace{3mm}

\noindent [R5] Ramakrishna, R., {\it Deforming an even representation II,
               raising the level}, Journal of Number Theory,
               72 (1998), 92-109.

\vspace{3mm} 

\noindent [S] Serre, J.-P., {Abelian and $l$-adic representations
                and elliptic curves}, W. A. Benjamin Inc., 1968.
\vspace{3mm}

\noindent [T1] Taylor, R., {\it Remarks on a conjecture of Fontaine and Mazur},
            Journal de l'Institut de Math. de Jussieu, 1 (2002). 1-19.

\vspace{3mm}

\noindent [T2] Taylor, R., {\it On icosahedral Artin representations II },
           to appear in the American Journal of Mathematics.

\vspace{3mm}

\noindent [Wa] Washington, L., {\it Galois cohomology}, in [FLT].

\vspace{3mm}

\noindent [W] Wiles, A., {\it Modular elliptic curves and 
Fermat's last theorem}, Ann. of
Math. (2) 141 (1995), no. 3, 443--551. 

\vspace{5mm}        

\noindent{\it Addresses of the authors:}

\vspace{3mm}

\noindent CK: School of Mathematics, TIFR, Homi Bhabha Road, 
Mumbai 400 005, INDIA. e-mail: shekhar@math.tifr.res.in;
\newline \noindent
155 S 1400 E, Dept of Math, Univ of Utah, Salt Lake City, UT 84112, USA:
e-mail: shekhar@math.utah.edu
\vspace{3mm}

\noindent RR: Department of Mathematics, Cornell University, Malott Hall,
Ithaca, NY 14853, USA. e-mail: ravi@math.cornell.edu
\end{document}